\documentclass{amsart}
\usepackage{amssymb}
\usepackage[all]{xy}
\usepackage{amsmath}
\usepackage{amsthm}

\theoremstyle{plain}
\newtheorem{lemma}{Lemma}[section]
\newtheorem{proposition}[lemma]{Proposition}
\newtheorem{theorem}[lemma]{Theorem}
\newtheorem{corollary}[lemma]{Corollary}
\theoremstyle{definition}
\newtheorem{definition}[lemma]{Definition}
\theoremstyle{remark}
\newtheorem{example}[lemma]{Example}
\newtheorem{remark}[lemma]{Remark}

\newtheorem{convention}[lemma]{Convention}

\newcommand{\longto}[1][]{\stackrel{#1}{\longrightarrow}}
\newcommand{\frakM}{\mathfrak{M}}
\newcommand{\calM}{\mathcal{M}}
\newcommand{\calU}{\mathcal{U}}
\newcommand{\calW}{\mathcal{W}}
\newcommand{\calE}{\mathcal{E}}
\newcommand{\calL}{\mathcal{L}}
\newcommand{\calC}{\mathcal{C}}
\newcommand{\calO}{\mathcal{O}}
\newcommand{\rmH}{\mathrm{H}}
\newcommand{\et}{\mathrm{\acute{e}t}}
\newcommand{\Het}{\rmH_{\et}}
\newcommand{\mods}[1]{{}^{\Z}/_{#1}}
\newcommand{\stab}{\mathrm{s}}
\newcommand{\semistab}{\mathrm{ss}}
\newcommand{\regstab}{\mathrm{rs}}
\newcommand{\pr}{\mathrm{pr}}
\newcommand{\id}{\mathrm{id}}
\newcommand{\ev}{\mathrm{ev}}
\newcommand{\coroots}{\mathrm{coroots}}

\newcommand{\dbar}{\bar{d}}
\newcommand{\Gbar}{\bar{G}}
\newcommand{\Gad}{G^{\ad}}
\newcommand{\gLie}{\mathfrak{g}}

\newcommand{\pLie}{\mathfrak{p}}
\newcommand{\Ga}{\mathbb{G}_a}
\newcommand{\Gm}{\mathbb{G}_{\mathrm{m}}}
\newcommand{\Z}{\mathbb{Z}}
\newcommand{\Q}{\mathbb{Q}}
\newcommand{\Frob}{\mathrm{Frob}}

\DeclareMathOperator{\Pic}{Pic}

\DeclareMathOperator{\Spec}{Spec}
\DeclareMathOperator{\Aut}{Aut}
\DeclareMathOperator{\End}{End}
\DeclareMathOperator{\Hom}{Hom}
\DeclareMathOperator{\NS}{NS}
\DeclareMathOperator{\GL}{GL}
\DeclareMathOperator{\PGL}{PGL}
\DeclareMathOperator{\ad}{ad}
\DeclareMathOperator{\weight}{wt}
\DeclareMathOperator{\rank}{rank}
\DeclareMathOperator{\coker}{coker}
\DeclareMathOperator{\characteristic}{char}

\begin{document}
\bibliographystyle{plain}

\title{Poincar\'{e} families of $G$-bundles on a curve}

\author[I. Biswas]{Indranil Biswas}
\address{School of Mathematics, Tata Institute of Fundamental Research, Homi Bhabha Road, Mumbai 400005, India}
\email{indranil@math.tifr.res.in}

\author[N. Hoffmann]{Norbert Hoffmann}
\address{Mathematisches Institut der Freien Universit\"at, Arnimallee 3, 14195 Berlin, Germany}
\email{norbert.hoffmann@fu-berlin.de}


\date{}

\begin{abstract}
Let $G$ be a reductive group over an algebraically closed field $k$.
Consider the moduli space of stable principal $G$-bundles on a smooth projective curve $C$ over $k$.
We give necessary and sufficient conditions for the existence of Poincar\'e bundles over open subsets of this moduli space,
and compute the orders of the corresponding obstruction classes.
This generalizes the previous results of Newstead, Ramanan and Balaji-Biswas-Nagaraj-Newstead to all reductive groups,
to all topological types of bundles, and also to all characteristics.
\end{abstract}

\thanks{The first author thanks the Freie Universit\"at Berlin for hospitality. 
The second author was supported by the SFB 647: Raum - Zeit - Materie.}

\maketitle

\section{Introduction}

Let $C$ be a compact Riemann surface of genus $g \geq 2$. Let $\frakM_{n, d}$ denote the coarse moduli space
of stable vector bundles $E$ over $C$ of fixed rank $n \geq 2$ and degree $d \in \Z$. Newstead proved
that there is no Poincar\'e bundle over $C \times \frakM_{n, 0}$, more precisely, the universal projective bundle
over $C \times \frakM_{n, 0}$ cannot be lifted to a topological vector bundle \cite{Ne}.
Ramanan proved that if $n$ and $d$ have a common divisor, then there is no Poincar\'e bundle over $C \times U$
for any Zariski-open subset $\emptyset \neq U \subseteq \frakM_{n, d}$ \cite{Ra}. If $n$ is coprime to $d$,
it is easy to see that there is a Poincar\'e bundle over $C \times \frakM_{n, d}$.

Now let $G$ be a complex semisimple group with one simple factor. Let $\frakM^{0, \regstab}_G$ be the moduli space of topologically trivial
regularly stable principal $G$-bundles over $C$ (a stable $G$-bundle is called regularly stable if its automorphism group coincides with the center of $G$).
The main theorem of \cite{BBNN} says that there is a Poincar\'{e} bundle (also called universal principal $G$-bundle) over $C \times \frakM^{0, \regstab}_G$
if and only if the center of $G$ is trivial. The proof in \cite{BBNN} is a generalization of the proof in \cite{Ne}.

In this paper, we address the existence of such Poincar\'{e} bundles for all topological types, not just topologically trivial.
The definition of Poincar\'{e} bundle (or Poincar\'{e} family) is recalled in Definition \ref{def:poinc}.
We give necessary and sufficient conditions for their existence over the whole regularly stable locus (Corollary \ref{cor:ob_generic}),
as well as over arbitrarily small Zariski-open subsets of it (Corollary \ref{cor:ob_global}).
In fact we compute the orders of the corresponding obstruction classes, which live in appropriate Brauer groups;
cf. Theorem \ref{thm:ob_generic} and Theorem \ref{thm:ob_global}.
These orders are given in terms of the root system of $G$ (see Section \ref{sec:obstructions}).
They are easy to compute for any given reductive group $G$.
The results are valid also in positive characteristics.

The proof of our results is based on the observation that the regularly stable locus in the moduli stack of principal $G$-bundles is a gerbe over the corresponding
locus in the coarse moduli scheme. The Brauer class of this gerbe is precisely the obstruction class mentioned above.
We have studied the Picard groups of such moduli stacks in \cite{pic}. This is used in Section \ref{sec:gerbes} to determine the order of the obstruction class.

In order to use this method, we need to know that the regularly stable locus is at least non-empty. This is already known in characteristic zero,
but apparently not in positive characteristic. We recall the notion of regular stability in Section \ref{sec:regstab},
and formulate the precise statement that we need as Theorem \ref{thm:simple}.
Its proof is given in Section \ref{sec:proof}, using some preparations in Sections \ref{sec:Brill-Noether} and \ref{sec:B-modules}.

\section{Regularly stable $G$-bundles} \label{sec:regstab}
Let $C$ be a connected smooth projective algebraic curve of genus $g \geq 2$ over an algebraically closed base field $k$.
Throughout this text, we assume that $G$ is a smooth connected linear algebraic group over $k$. The Lie algebra of $G$
will be denoted by $\gLie$.

Let $\calM_G$ denote the moduli stack of principal $G$-bundles on $C$. One knows that $\calM_G$ is an Artin stack, locally of finite type and smooth over $k$.
Its decomposition into connected components
\begin{equation*}
  \calM_G = \coprod_{d \in \pi_1( G)} \calM_G^d
\end{equation*}
is indexed by the fundamental group $\pi_1( G)$ of $G$; see \cite[Section 5]{pi0}.
Let $G_u \subseteq G$ be the unipotent radical. Choosing a maximal torus $T$ in the reductive quotient $G/G_u$, the group $\pi_1( G)$ is by definition
the quotient of $\Hom( \Gm, T)$ modulo its subgroup generated by the coroots of $G/G_u$.

Assume that $G$ is reductive. According to \cite[Proposition 3.20]{yogish}, the stable and the semistable principal $G$-bundles form open substacks
\begin{equation*}
  \calM_G^{d, \stab} \subseteq \calM_G^{d, \semistab} \subseteq \calM_G^d.
\end{equation*}
\begin{lemma} \label{generic_stable}
  The complement of $\calM_G^{d, \stab}$ has codimension $\geq g-1$ in $\calM_G^d$.
\end{lemma}
\begin{proof}
  Let $\calC$ be an irreducible component of the complement $\calM_G^d \setminus \calM_G^{d, \stab}$. Let $K$ be an algebraically closed field
  containing $k$, and let $E$ be a principal $G$-bundle on $C \otimes_k K$ whose classifying morphism $\Spec( K) \to \calM_G$ maps onto the generic
  point of $\calC$. Then $E$ is not stable, so there is a parabolic subgroup $P \subseteq G$ and a reduction of $E$ to a principal $P$-bundle $F$ with
  $\deg( \ad( F)) \leq 0$; here $\ad( F) := F \times^P \pLie$ denotes the vector bundle obtained from $F$ via the adjoint representation
  $\ad: P \to \Aut( \pLie)$, where $\pLie$ is the Lie algebra of $P$. We consider the connected component
  \begin{equation*}
    \calM_P^{\delta} \subseteq \calM_P
  \end{equation*}
  that contains $F$. Extending the structure group defines a representable $1$-morphism $\calM_P^{\delta} \to \calM_G^d$ 
  which is dominant onto $\calC$. Hence
  \begin{equation*}
    \dim \calC \leq \dim \calM_P^{\delta} = (g-1) \dim P - \deg( \ad( F)) \leq (g-1) \dim P.
  \end{equation*}
  So the codimension of $\calC$ in $\calM_G^d$ is $\geq (g-1) (\dim G - \dim P) \geq g-1$.
\end{proof}
Due to \cite{GLSS,GLSS2} and \cite{jochen,jochen2}, there is a projective coarse moduli scheme
\begin{equation*}
  \pi: \calM_G^{d, \semistab} \longto \frakM_G^d
\end{equation*}
of semistable principal $G$-bundles $E$ on $C$ of type $d \in \pi_1( G)$. These constructions also show that
the stable locus $\calM_G^{d, \stab} \subseteq \calM_G^{d, \semistab}$ is the inverse image of an open subscheme $\frakM_G^{d, \stab} \subseteq \frakM_G^d$.

Let $Z \subseteq G$ denote the (scheme-theoretic) center. It is a group scheme of the form $Z \cong \Gm^r \times \mu_{n_1} \times \cdots \times \mu_{n_s}$
with $r, s \geq 0$ and $n_1, \ldots, n_s \geq 1$.
\begin{definition} \label{regstab}
  A stable principal $G$-bundle $E$ on $C$ is \emph{regularly stable} if the canonical morphism
  \begin{equation*}
    Z \longto \Aut( E)
  \end{equation*}
  into the group scheme of global automorphisms of $E$ is an isomorphism.
\end{definition}
\begin{proposition} \label{prop:regstab}
  The locus of regularly stable principal $G$-bundles
  \begin{equation*}
    \calM_G^{d, \regstab} \subseteq \calM_G^{d, \stab}
  \end{equation*}
  is an open substack. It is the inverse image of an open subscheme
  \begin{equation*}
    \frakM_G^{d, \regstab} \subseteq \frakM_G^{d, \stab}.
  \end{equation*}
\end{proposition}

In order to prove this proposition, we recall Luna's \'{e}tale slice theorem \cite{luna},
and its generalization to positive characteristic in \cite{bardsley-richardson}.

Let $H$ be a smooth reductive algebraic group over $k$. Suppose that $H$ acts on an affine variety $X$ of finite type over $k$. 
Let $X/H$ denote the GIT-quotient. Let $x \in X( k)$ be a closed point. Its orbit $H \cdot x$ is a smooth subvariety of $X$.
The orbit map $H \twoheadrightarrow H \cdot x$ is an fppf-locally trivial principal bundle under the scheme-theoretic stabilizer $H_x \subseteq H$.

An $H_x$-stable affine subvariety $S \subseteq X$ with $x \in S$ is called an \emph{\'{e}tale slice} at $x$ if the commutative diagram of canonical maps
\begin{equation*} \xymatrix{
  H \times^{H_x} S \ar[r] \ar[d] & X \ar[d]\\
  (H \times^{H_x} S)/H \cong S/H_x \ar[r] & X/H 
} \end{equation*}
is cartesian, and the horizontal maps are both \'{e}tale.
\begin{proposition} \label{prop:slice}
  If the orbit $H \cdot x$ is closed in $X$ and the stabilizer $H_x$ is linearly reductive, then there exists an \'{e}tale slice at $x$.
\end{proposition}
\begin{proof}
  Bardsley and Richardson have proved this under the assumption that the orbit map $H \twoheadrightarrow H \cdot x$ is separable
  \cite[Proposition 7.6]{bardsley-richardson}. However, this assumption is not needed in their argument if one works with the (possibly non-reduced)
  scheme-theoretic stabilizer $H_x$.

  More precisely, \cite[Proposition 7.5]{bardsley-richardson} allows us to assume without loss of generality that $X$ is a vector space over $k$,
  and that $H$ acts linearly on $X$. Since $H_x$ is linearly reductive, there is an $H_x$-stable affine linear subspace $S \subseteq X$ containing $x$
  such that $S$ and the orbit $H \cdot x$ are transversal at $x$. This means that the canonical morphism of smooth varieties
  \begin{equation*}
    H \times^{H_x} S \longto X
  \end{equation*}
  is \'{e}tale at the point $(1, x) \in H \times^{H_x} S$. Using \cite[Theorem 6.2]{bardsley-richardson}, this implies that some Zariski-open
  neighborhood of $x$ in $S$ is an \'{e}tale slice at $x$.
\end{proof}
\begin{proof}[Proof of Proposition \ref{prop:regstab}]
  According to the construction in \cite{GLSS,GLSS2}, the stable locus $\frakM_G^{d, \stab}$ is a GIT-quotient of a quasiprojective scheme $X$ over $k$
  modulo a reductive group $H$. More precisely, the quotient morphism
  \begin{equation*}
    X \longto \frakM_G^{d, \stab}
  \end{equation*}
  is affine, its fibers are precisely the $H$-orbits, the $H$-invariant open subschemes of $X$ are precisely the inverse images of
  open subschemes in $\frakM_G^{d, \stab}$, and the stack quotient $[X/H]$ is the moduli stack $\calM_G^{d, \stab}$.

  In particular, every orbit $H \cdot x \subseteq X$ with $x \in X( k)$ is closed, and the scheme-theoretic stabilizer $H_x$ always contains $Z$.
  If $H_x$ is just $Z$, then it is in particular linearly reductive \cite[II, \S2, 2.5]{demazure-gabriel}, so Proposition \ref{prop:slice} applies
  and shows that there is an \'{e}tale slice at $x$.

  It follows that the set of all points $x \in X( k)$ with stabilizer $H_x = Z$ is Zariski-open.
  Since it is also $H$-invariant, it defines an open substack $\calM_G^{d, \regstab} \subseteq \calM_G^{d, \stab}$ and an open subscheme
  $\frakM_G^{d, \regstab} \subseteq \frakM_G^{d, \stab}$.
\end{proof}
\begin{theorem} \label{thm:simple}
  We still assume that $G$ is reductive.
  \begin{itemize}
   \item[i)] If $g \geq 3$, then $\calM_G^{d, \regstab}$ is non-empty.
   \item[ii)] If $g \geq 4$, then the complement of $\calM_G^{d, \regstab}$ in $\calM_G^d$ has codimension at least two.
  \end{itemize}
\end{theorem}

In characteristic $0$, this is proved in \cite[Theorem II.6]{faltings}, with better bounds on $g$.
However, we were not able to generalize the cited proof to positive characteristic, mainly due to difficulties with finite unipotent automorphism groups.
After some preparations in Sections \ref{sec:Brill-Noether} and \ref{sec:B-modules}, we prove Theorem \ref{thm:simple} in Section \ref{sec:proof}.

\section{Brill-Noether loci} \label{sec:Brill-Noether}
By a $G$-module $V$, we mean a finite-dimensional vector space $V$ over $k$ together with a morphism $G \to \GL( V)$ of algebraic groups over $k$. 
Given such a $G$-module $V$, we can associate to each principal $G$-bundle $E$ on $C$ the vector bundle $E \times^G V$ on $C$ with fiber $V$.
The Brill-Noether locus
\begin{equation*}
  \calW_V := \{[E]: \rmH^0( C, E \times^G V) \neq 0\} \subseteq \calM_G
\end{equation*}
is closed, according to the semicontinuity theorem. This section deals with codimension estimates for such loci.

Each group homomorphism $\varphi: H \to G$ induces, by extension of the structure group, a $1$-morphism $\varphi_*: \calM_H \to \calM_G$.
Moreover, $\varphi$ allows us to turn $G$-modules $V$ into $H$-modules.
We note that $\calW_V \subseteq \calM_H$ is the inverse image of $\calW_V \subseteq \calM_G$ in this situation.

For a character $\chi: G \to \Gm$, let $\ell_{\chi}$ be the resulting $1$-dimensional $G$-module.
Let $\langle \chi, d \rangle \in \Z$ denote the image of $d \in \pi_1( G)$ under the group homomorphism $\chi_*: \pi_1( G) \to \pi_1( \Gm) = \Z$.
\begin{lemma} \label{lemma:chi}
  If $\chi \neq 0$, then $\calW_{\ell_{\chi}}$ has codimension $\geq g - \langle \chi, d \rangle$ in $\calM_G^d$.
\end{lemma}
\begin{proof}
  We start with the special case $G = \Gm$ and $\chi = \id$. Note that $\calW_{\ell_{\id}} \subseteq \calM_{\Gm}^d$ is empty for $d < 0$,
  and  for $d \geq 0$ it corresponds to the image of the Abel-Jacobi map $C^d \to \Pic^d( C)$.
  This image has dimension $\leq d$, and hence codimension $\geq g - d$ in $\Pic^d( C)$.
  Thus $\calW_{\ell_{\id}}$ has codimension $\geq g - d$ in $\calM_{\Gm}^d$, as claimed.

  In the general case, $\calW_{\ell_{\chi}} \subseteq \calM_G$ is the inverse image of $\calW_{\ell_{\id}} \subseteq \calM_{\Gm}$
  under the induced $1$-morphism $\chi_*: \calM_G \to \calM_{\Gm}$. Using the above special case, it thus suffices to show that this $1$-morphism $\chi_*$ is flat.

  If a closed normal subgroup $N \subseteq G$ is smooth, then the canonical $1$-morphism $\calM_G \to \calM_{G/N}$ is also smooth by deformation theory.
  We apply this to the reduced identity component $N := \ker( \chi)^0 \subseteq G$.
  It contains the unipotent radical $G_u \subseteq G$ and the commutator subgroup $[G, G] \subseteq G$, so $G/N$ is a torus.
  The character $G/N \to \Gm$ induced by $\chi$ is nontrivial and has finite kernel, so it is an isogeny.
  In particular, $G/N$ has rank $1$, and the induced $1$-morphism $\calM_{G/N} \to \calM_{\Gm}$ is flat.
  It follows that the composition $\chi_*: \calM_G \to \calM_{G/N} \to \calM_{\Gm}$ is also flat, as required.
\end{proof}
For the rest of this section, we consider the group $G := \Ga \rtimes \Gm$ for the standard action of $\Gm$ on $\Ga$.
Let $\pi: G \twoheadrightarrow \Gm$ denote the projection.

Put $q = 1$, or let $q$ be a power of $p$ if $\characteristic( k) = p \geq 2$.
Let $\wp^q$ denote the $2$-dimensional $G$-module given by the homomorphism $G \to \GL_2$ that sends $t \in \Ga$ to
$\big( \begin{smallmatrix} 1 & 0\\t^q & 1 \end{smallmatrix} \big)$, and $h \in \Gm$ to $\big( \begin{smallmatrix} 1 & 0\\0 & h^q \end{smallmatrix} \big)$.
We have an exact sequence
\begin{equation*}
  0 \longto \ell_{q \pi} \longto \wp^q \longto \ell_0 \longto 0
\end{equation*}
of $G$-modules. In the case $q = 1$, we write $\wp := \wp^1$; here we get in particular an exact sequence of $G$-modules
\begin{equation*}
  0 \longto \ell_{\pi} \longto \wp \longto \ell_0 \longto 0.
\end{equation*}
It yields an equivalence of categories between principal $G$-bundles on $C$ and exact sequences of vector bundles on $C$
\begin{equation*}
  0 \longto L \longto E \longto \calO_C \longto 0
\end{equation*}
with $\rank( L) = 1$ and $\rank( E) = 2$.
Thus the fiber of $\pi_*: \calM_G \to \calM_{\Gm}$ over a line bundle $L$ on $C$ parameterizes extensions of $\calO_C$ by $L$.

The type $d \in \pi_1( G) = \Z$ of a principal $G$-bundle on $C$ is the degree of the corresponding line bundle $L$.
\begin{lemma} \label{lemma:b}
  The locus $\calW_{\wp} \setminus \calW_{\ell_{\pi}} \subseteq \calM_G$ has dimension $g - 1$.
\end{lemma}
\begin{proof}
  Let a line bundle $L$ on $C$ with $\rmH^0( C, L) = 0$ be given, and an extension $E$ of $\calO_C$ by $L$ as above.
  Then $\rmH^0( C, E) \neq 0$ if and only if this extension splits.
  Hence the restriction of the $1$-morphism $\pi_*$ is an isomorphism of $\calW_{\wp} \setminus \calW_{\ell_{\pi}}$ onto $\calM_{\Gm} \setminus \calW_{\ell_{\id}}$.
\end{proof}
\begin{corollary} \label{cor:b}
  If $d \leq 0$, then $\calW_{\wp}$ has codimension $\geq g-1$ in $\calM^d_G$.
\end{corollary}
\begin{proof}
  Since its fibers parameterize extensions by $\calO_C$, the $1$-morphism $\pi_*: \calM_G^d \to \calM_{\Gm}^d$ is smooth of relative dimension $g - 1 - d$;
  see Lemma 2.10 of \cite{par}. Hence $\calM_G^d$ has dimension $2g - 2 - d$.

  If $d \leq 0$, then $\calW_{\wp} \setminus \calW_{\ell_{\pi}}$ has codimension $\geq g-1$ in $\calM^d_G$ by Lemma \ref{lemma:b}, and
  $\calW_{\ell_{\pi}}$ has codimension $\geq g$ in $\calM_G^d$ by Lemma \ref{lemma:chi}.
\end{proof}
By a \emph{prime divisor} $D$ in an integral Artin stack $\calM$, we mean an integral closed substack $D \subseteq \calM$ of codimension $1$.
\begin{lemma} \label{lemma:Wwp}
  Suppose $\characteristic( k) = p \geq 2$, and $d \leq 0$.
  \begin{itemize}
   \item[i)] If $g \geq 2$, then the complement of $\calW_{\wp^p}$ in $\calM_G^d$ is non-empty.
   \item[ii)] If $g \geq 3$, and some prime divisor $D \subseteq \calM_G^d$ is contained in $\calW_{\wp^p}$,
    then $D = (\pi_*)^{-1}( D')$ for some prime divisor $D' \subseteq \calM_{\Gm}^d$.
   \item[iii)] If $p = 2$, then $\calW_{\wp^2}$ has codimension $\geq g-1$ in $\calM_G^d$.
  \end{itemize}
\end{lemma}
\begin{proof}
  We consider a principal $G$-bundle on $C$ whose moduli point is in $\calW_{\wp^p}$, but not in $\calW_{\wp} \cup \calW_{\ell_{p \pi}}$.
  The exact sequences $0 \to \ell_{\pi} \to \wp \to \ell_0 \to 0$ and $0 \to \ell_{p \pi} \to \wp^p \to \ell_0 \to 0$ of $G$-modules
  induce exact sequences
  \begin{equation*}
    0 \to L \longto E \longto \calO_C \to 0 \quad\text{and}\quad 0 \to L^{\otimes p} \longto \Frob^*( E) \longto \calO_C \to 0
  \end{equation*}
  of associated vector bundles on $C$, where $\Frob: C \to C$ is the absolute Frobenius.
  Our assumption means that the extension $E$ doesn't split, but the extension $\Frob^*( E)$ does split.

  In particular, there is a nonzero morphism $\Frob^*( E) \to L^{\otimes p}$ of vector bundles on $C$.
  Using adjunction, we obtain a nonzero morphism
  \begin{equation*}
    \varphi: E \longto \Frob_*( L^{\otimes p}) = \Frob_*( \Frob^*( L)) = L \otimes \Frob_*( \calO_C).
  \end{equation*}  
  Since $\rmH^0( C, \Frob_*( \calO_C)) = \rmH^0( C, \calO_C) = k$, the restriction of $\varphi$ maps $L \subseteq E$ to the subbundle 
  $L \subseteq L \otimes \Frob_*( \calO_C)$. Thus we obtain a nonzero morphism of short exact sequences over $C$
  \begin{equation} \label{eq:ext_map} \xymatrix{
    0 \ar[r] & L \ar[r] \ar[d] & E \ar[r] \ar[d]^{\varphi} & \calO_C \ar[r] \ar[d]^{\overline{\varphi}} & 0\\
    0 \ar[r] & L \ar[r] & L \otimes \Frob_*( \calO_C) \ar[r] & L \otimes \Frob_*( \calO_C)/\calO_C \ar[r] & 0.
  } \end{equation}
  The top row doesn't split, so $\overline{ \varphi} \neq 0$. This shows
  \begin{equation*}
    \calW_{\wp^p} \subseteq \calW_{\wp} \cup \calW_{\ell_{p \pi}} \cup ( \pi_*)^{-1}( \calW')
  \end{equation*}
  as closed loci in $\calM_G$, where we put
  \begin{equation*} \label{eq:W'}
    \calW' := \{[ L]: \rmH^0 \big( C, L \otimes \Frob_*( \calO_C)/\calO_C \big) \neq 0\} \subseteq \calM_{\Gm}.
  \end{equation*}
  As $d\leq 0$ by assumption, $\calW_{\wp} \cup \calW_{\ell_{p \pi}}$ has codimension $\geq g-1$ in $\calM_G^d$ due to Corollary \ref{cor:b} and Lemma \ref{lemma:chi}.
  Due to \cite[Th\'{e}or\`{e}me 4.1.1]{raynaud}, we have $\calM_{\Gm}^0 \not\subseteq \calW'$ and hence $\calM_{\Gm}^d \not\subseteq \calW'$.
  This implies part (i) of the lemma.

  In the situation of part (ii), we conclude $D \subseteq (\pi_*)^{-1}( D')$ for some irreducible component $D'$ of $\calW' \cap \calM_{\Gm}^d$.
  But $D$ is a prime divisor, $\pi_*$ is smooth, and $D' \neq \calM_{\Gm}^d$. It follows that $D'$ is also a prime divisor, and that $D = (\pi_*)^{-1}( D')$.
  This proves part (ii) of the lemma.

  Now suppose $p = 2$. Then $\Frob_*( \calO_C)/\calO_C$ is a line bundle on $C$, of degree $g - 1$ \cite[p. 118]{raynaud}.
  If $\deg( L) = d$, then $\overline{ \varphi}$ is a nonzero section of a line bundle of degree $g - 1 + d$.
  Sending our $G$-bundle to the divisor of this section defines a morphism
  \begin{equation*}
    \calM_G^d \cap \calW_{\wp^2} \setminus (\calW_{\wp} \cup \calW_{\ell_{2 \pi}}) \longto C^{(g-1+d)}.
  \end{equation*}
  This is an open embedding, since the top row in \eqref{eq:ext_map} can be reconstructed from $L$ and $\overline{\varphi}$
  by pulling back the bottom row. Hence we conclude
  \begin{equation*}
    \dim \big( \calM_G^d \cap \calW_{\wp^2} \setminus (\calW_{\wp} \cup \calW_{\ell_{2 \pi}}) \big) \leq g - 1 + d.
  \end{equation*}
  As we are assuming $d \leq 0$, part (iii) of the lemma follows. 
\end{proof}

\section{$B$-modules} \label{sec:B-modules}

In this section, we assume that our algebraic group $G$ over $k$ is of adjoint type and simple.
Choose a Borel subgroup $B \subseteq G$, and a maximal torus $T \subseteq B$. Let us denote by
\begin{equation*}
  \Hom( T, \Gm) \supset \Phi \supset \Phi^+ \ni \alpha_1, \ldots, \alpha_l
\end{equation*}
the root system of $G$ with respect to $T$, the positive roots with respect to $B$, and the simple roots, respectively. We adopt the following:
\begin{convention} \label{conv:G2}
  In the case $\Phi \cong G_2$, we number the two simple roots $\alpha_1, \alpha_2 \in \Phi^+$ in such a way that $\alpha_1$ is short and $\alpha_2$ is long.
\end{convention}
Every $B$-module $V$ decomposes into $T$-eigenspaces 
\begin{equation*}
  V = \bigoplus_{\chi \in \Hom( T, \Gm)} V_{\chi}.
\end{equation*}
E. g. the adjoint action of $B \subseteq G$ on $\gLie$ yields the Cartan decomposition
\begin{equation*}
  \gLie = \gLie_0 \oplus \bigoplus_{\alpha \in \Phi} \gLie_{\alpha}.
\end{equation*}
We choose basis vectors $e_{\alpha} \in \gLie_{\alpha}$ that form a Chevalley system; see for example \cite[p. 56f.]{carter}.
Given a root $\alpha \in \Phi$, \cite[Expos\'e XXII, Th\'{e}or\`{e}me 1.1]{sga3.3} provides a unique $T$-equivariant morphism of algebraic groups
\begin{equation*}
  \exp_{\alpha}: \gLie_{\alpha} \longto G
\end{equation*}
whose derivative at $0$ is the inclusion $\gLie_{\alpha} \hookrightarrow \gLie$.
The map $\exp_{\alpha}$ is an isomorphism onto a closed subgroup, which we denote by $U_{\alpha} \subseteq G$.

Put $u_{\alpha} := \exp_{\alpha}( e_{\alpha}) \in U_{\alpha}$.
Let $U$ denote the unipotent radical of $B$; it is generated by the $U_{\alpha}$ with $\alpha \in \Phi^+$.
For every $B$-module $V$, the action of $U_{\alpha} \subseteq B$ on $T$-eigenvectors $v \in V_{\chi}$ satisfies
\begin{equation} \label{eq:U_alpha}
  \exp_{\alpha}( t e_{\alpha}) \cdot v = v + \sum_{n \geq 1} t^n \rho_{\chi, n \alpha}( v)
\end{equation}
for all $t \in k$, with linear maps $\rho_{\chi, n \alpha}: V_{\chi} \to V_{\chi + n \alpha}$; see e.\,g. Lemma 5.2 of \cite{borel}.
Hence the space of invariants $V^B \subseteq V$ is the kernel of the linear map
\begin{equation*}
  \oplus \rho_{0, n \alpha}: V_0 \longto \bigoplus_{\alpha \in \Phi^+} \bigoplus_{n \geq 1} V_{n \alpha}.
\end{equation*}
Therefore, given a short exact sequence $0 \to V' \to V \to V'' \to 0$ of $B$-modules, the snake lemma yields an exact sequence of vector spaces
\begin{equation} \label{eq:VB}
  0 \to (V')^B \longto V^B \longto (V'')^B
    \longto \coker \big( V_0' \xrightarrow{\oplus \rho_{0, n \alpha}} \bigoplus_{\alpha \in \Phi^+} \bigoplus_{n \geq 1} V_{n \alpha}' \big).
\end{equation}

For each simple root $\alpha_i$, we define a homomorphism of algebraic groups
\begin{equation*}
  \pi_i: B = U \rtimes T \longto \Ga \rtimes \Gm
\end{equation*}
as the product of $\alpha_i: T \to \Gm$ with the projection $U \to \Ga$ that vanishes on all $U_{\alpha}$ with $\alpha \in \Phi^+ \setminus \{\alpha_i\}$,
and maps $\exp( t e_{\alpha_i}) \in U_{\alpha_i}$ to $t \in \Ga$.

As $G$ is of adjoint type, the product map $\prod_i \alpha_i: T \to \Gm^l$ is an isomorphism. Thus we get an exact sequence
\begin{equation} \label{eq:U'}
  1 \longto U' \longto B \xrightarrow{\prod_i \pi_i} ( \Ga \rtimes \Gm)^l \longto 1
\end{equation}
of algebraic groups over $k$, where $U'$ is generated by all $U_{\alpha}$ with $\alpha \in \Phi^+$ not simple.
In particular, $U'$ is smooth and connected.

Put $q = 1$, or let $q$ be a power of $p$ if $\characteristic( k) = p \geq 2$.
Let $\wp_i^q$ denote the $B$-module obtained via $\pi_i$ from the $(\Ga \rtimes \Gm)$-module $\wp^q$, and put $\wp_i := \wp_i^1$.
These are indecomposable $B$-modules of dimension $2$.
\begin{proposition} \label{prop:wp}
  Let $V$ be an indecomposable $B$-module of dimension $2$.
  Then $V \cong \ell_{\chi} \otimes \wp_i^q$ for some character $\chi: B \to \Gm$, some simple root $\alpha_i$, and some number $q$ as above.
\end{proposition}
\begin{proof}
  Because $V$ is indecomposable of dimension $2$, there is a positive root $\alpha$ such that $U_{\alpha} \subseteq B$ acts nontrivially on $V$.

  Using relation \eqref{eq:U_alpha}, this observation implies $V \cong \ell_{\chi} \oplus \ell_{\chi + q \alpha}$ as $T$-modules,
  for some character $\chi: T \to \Gm$ and some integer $q \geq 1$.

  Choosing eigenvectors $v_{\psi} \in V_{\psi}$, we moreover see that $\exp( t e_{\alpha}) \in U_{\alpha}$ maps $v_{\chi}$ to $v_{\chi} + c t^q v_{\chi + q \alpha}$
  for a constant $c \in k$. Here $c \neq 0$, as $V$ is indecomposable. Replacing $v_{\chi}$ by $c v_{\chi}$ yields $c = 1$ without loss of generality.

  Since the map $U_{\alpha} \to \GL( V)$ at hand is a group homomorphism, we have $q = 1$, or $q$ is a power of $\characteristic( k) \geq 2$.
  We claim that $\alpha$ is a simple root $\alpha_i$; then $V \cong \ell_{\chi} \otimes \wp_i^q$ will follow.

  Suppose that $\alpha$ is not simple. Then $\alpha = \beta + \gamma$ for some $\beta, \gamma \in \Phi^+$.
  We assume without loss of generality that $\gamma$ is not longer than $\beta$.
  Let $n \geq 1$ be maximal such that $\beta_0 := \alpha - n \gamma \in \Phi^+$; then $n \leq 3$.

  If $\gamma$ is long, then $\beta$ is also long by assumption, so $n = 1$, and $\beta_0 = \beta$ is long as well.
  Thus $\gamma$ is not longer than $\beta_0$ in any case.

  Due to \cite[Expos\'e XXIII, Proposition 6.4]{sga3.3}, the equation
  \begin{equation*}
    u_{\beta_0} u_{\gamma} = u_{\gamma} u_{\beta_0} (u_{\beta_0 + \gamma})^{\pm 1} (u_{\beta_0 + 2\gamma})^{\pm 1} (u_{\beta_0 + 3\gamma})^{\pm 1} (u_{2\beta_0 + 3\gamma})^{\pm1}
  \end{equation*}
  holds in $B$, using the convention $u_{\psi} := 1$ for $\psi \not\in \Phi$.
  One of these factors is $(u_{\alpha})^{\pm 1}$; it acts nontrivially on $V$ by assumption.
  But all other factors in this equation act trivially on $V$, due to relation \eqref{eq:U_alpha}.

  This contradiction proves the claim.
\end{proof}
\begin{corollary} \label{cor:VB0}
  If $V$ is a $B$-module with $V^B = 0$,
  then $V$ is a successive extension of $B$-modules of the form $\ell_{\chi}$ with $\chi \neq 0$, or of the form $\wp^q_i$.
\end{corollary}
\begin{proof}
  Suppose $V^B = 0$ and $V \neq 0$. We argue by induction on $\dim V$. As $B$ is solvable, there is a $1$-dimensional $B$-submodule $V^1 \subseteq V$.
  From the exact sequence \eqref{eq:VB}, we see that there are two possible cases:
  \begin{itemize}
   \item $(V/V^1)^B = 0$, and $V^1 \cong \ell_{\chi}$ for some nontrivial character $\chi: T \to \Gm$.
   \item $\dim (V/V^1)^B = 1$, and $V^1 \cong \ell_{n \alpha}$ for some $\alpha \in \Phi^+$ and some $n \geq 1$.
  \end{itemize}
  In the second case, let $V^2 \subseteq V$ be the inverse image of $(V/V^1)^B \subseteq V/V^1$.
  Then $(V^2)^B = 0$, and $V^2 \cong \ell_0 \oplus \ell_{n \alpha}$ as $T$-modules.
  Due to Proposition \ref{prop:wp}, this implies $V^2 \cong \wp^q_i$ for some $i$ and $q$.
  The short exact sequence
  \begin{equation*}
    0 \longto (V/V^1)^B \longto V/V^1 \longto V/V^2 \longto 0
  \end{equation*}
  of $B$-modules yields an exact sequence \eqref{eq:VB} of vector spaces, from which we conclude $(V/V^2)^B = 0$.
  This completes the induction.
\end{proof}
\begin{proposition} \label{prop:EndB}
  Let $\phi: \gLie \to \gLie$ be a $B$-equivariant $k$-linear map. Then $\phi = \lambda \cdot \id_{\gLie}$ for some constant $\lambda \in k$.
\end{proposition}
\begin{proof}
  Since $\phi$ is a $T$-module endomorphism, it is a direct sum of components $\phi_0 \in \End( \gLie_0)$ and
  $\phi_{\alpha} \in \End(\gLie_{\alpha})$ for $\alpha \in \Phi$.
  As $\dim \gLie_{\alpha} = 1$, we have $\phi_{\alpha} = \lambda_{\alpha} \cdot \id_{\gLie_{\alpha}}$ with $\lambda_{\alpha} \in k$. 
  We put $\lambda := \lambda_{\theta}$ for the highest root $\theta \in \Phi^+$.

  Suppose that $\alpha = \beta + \gamma$ holds for $\alpha, \beta, \gamma \in \Phi^+$,
  and let $n \geq 1$ be maximal such that $\beta_0 := \alpha - n \gamma \in \Phi^+$. The action of $U_{\gamma}$ on $\gLie$ satisfies
  \begin{equation*}
    u_{\gamma} \cdot e_{\beta_0} = e_{\beta_0} + \sum_{m \geq 1} \pm e_{\beta_0 + m \gamma}
  \end{equation*}
  according to \cite[p. 64]{carter}, with the convention $e_{\psi} = 0$ for $\psi \not\in \Phi$.
  Comparing eigenvalues of $\phi$ on both sides, we conclude $\lambda_{\beta_0} = \lambda_{\beta_0 + \gamma} = \ldots$,
  and in particular $\lambda_{\beta} = \lambda_{\alpha}$. By symmetry, we also have $\lambda_{\gamma} = \lambda_{\alpha}$. 

  Writing $\alpha$ as a sum of simple roots $\alpha_{i_r}$, an iteration of this argument shows $\lambda_{\alpha_{i_r}} = \lambda_{\alpha}$ for all $r$.
  All simple roots appear in the highest root, so we get in particular $\lambda_{\alpha_i} = \lambda $ for all $i$,
  and hence $\lambda_{\alpha} = \lambda$ for all $\alpha \in \Phi^+$.

  For every $\alpha \in \Phi^+$, the formulas in \cite[p. 64]{carter} contain in particular
  \begin{equation*}
    u_{\alpha} \cdot e_{-\alpha} = e_{-\alpha} + [ e_{\alpha}, e_{-\alpha}] - e_{\alpha}.
  \end{equation*}
  Comparing the effect of $\phi$ on both sides shows $\lambda_{-\alpha} = \lambda_{\alpha} = \lambda$.
  
  The map $\gLie_0 \to \oplus_i \gLie_{\alpha_i}$, $h \mapsto (u_{\alpha_i} \cdot h - h)_i$, is an isomorphism due to the exact sequence \eqref{eq:U'}.
  Since $\phi$ is $B$-equivariant, it commutes with this isomorphism. 
  Hence $\phi_0 = \lambda \cdot \id_{\gLie_0}$ follows; thus $\phi = \lambda \cdot \id_{\gLie}$.
\end{proof}
\begin{corollary} \label{cor:End=ext}
  The $B$-module $V := \End( \gLie)/k \cdot \id$ is a successive extension of $B$-modules of the form
  \begin{itemize}
   \item $\ell_{\alpha}$ for some root $\alpha \in \Phi$, or
   \item $\ell_{\alpha - \beta}$ for some pair of different roots $\alpha \neq \beta \in \Phi$, or
   \item $\wp_i$ for some simple root $\alpha_i$, or
   \item $\wp^2_i$ for some simple root $\alpha_i$ in the case $\characteristic( k) = 2$, or
   \item $\wp^3_1$  in the case $\characteristic( k) = 3$ and $\Phi \cong G_2$.
  \end{itemize}
\end{corollary}
\begin{proof}
  We use the short exact sequence of $B$-modules
  \begin{equation*}
    0 \longto \ell_0 \cong k \cdot \id \longto \End( \gLie) \longto V \longto 0.
  \end{equation*}
  Proposition \ref{prop:EndB} states $\End( \gLie)^B = k \cdot \id$. Therefore, the associated exact sequence \eqref{eq:VB} of vector spaces shows $V^B=0$.
  Hence Corollary \ref{cor:VB0} applies, so $V$ is a successive extension of some $\wp^q_i$ and some $\ell_{\chi}$ for nontrivial $\chi$.

  Here $\ell_{\chi}$ can only appear if the $T$-eigenspace $V_{\chi}$ is nonzero, or in other words if $\chi$ is a root or a difference of roots.
  Similarly, $\wp^q_i$ can only appear if $q \alpha_i$ is a root or a difference of roots. According to the classification of root systems,
  the latter happens only if $q \leq 2$, or if $q = 3$ and $\Phi \cong G_2$ and $\alpha_i = \alpha_1$ is short.
\end{proof}

\section{Proof of Theorem \ref{thm:simple}} \label{sec:proof}
Let $G$ be reductive, with center $Z$, and let $\pi: G \to G/Z$ be the projection.
For each principal $G$-bundle $E$ with induced principal $(G/Z)$-bundle $\pi_* E$, we have $\Aut( E)/Z \subseteq \Aut( \pi_* E)$.
Thus $\calM_G^{\regstab}$ contains the inverse image of $\calM_{G/Z}^{\regstab}$ under the $1$-morphism $\pi_*: \calM_G \to \calM_{G/Z}$.

As this $\pi_*$ is flat according to \cite[Lemma 2.2.2]{pic}, it suffices to prove the theorem for $G/Z$ instead of $G$.
Thus we can assume without loss of generality that $G$ is of adjoint type, and also that $G$ is simple.

Let $E$ be a principal $G$-bundle on $C$. If its adjoint vector bundle $\ad( E) = E \times^G \gLie$ has only scalar endomorphisms, then $\Aut( E)$ is trivial.
This shows that the complement of $\calM_G^{\regstab}$ in $\calM_G^{\stab}$ is contained in the Brill-Noether locus
$\calW_{\End( \gLie)/k \cdot \id} \subseteq \calM_G$. Hence Theorem \ref{thm:simple} is a consequence of the following proposition.
\begin{proposition} \label{prop:W_End}
  Let $G$ be simple of adjoint type. Let $d \in \pi_1( G)$ be given.
  \begin{itemize}
   \item[i)] If $g \geq 3$, then the complement of $\calW_{\End( \gLie)/k \cdot \id}$ in $\calM_G^d$ is non-empty.
   \item[ii)] If $g \geq 4$, then $\calW_{\End( \gLie)/k \cdot \id}$ has codimension $\geq 2$ in $\calM_G^d$.
  \end{itemize}
\end{proposition}
In order to prove this proposition, we use the notation of the previous section. So $T \subseteq B \subseteq G$ is a maximal torus in a Borel subgroup, and
\begin{equation*}
  \alpha_1, \ldots, \alpha_l \in \Phi^+ \subseteq \Phi \subseteq \Hom( T, \Gm)
\end{equation*}
are the simple roots, the positive roots, and all roots of $G$, respectively.
If $G$ is of type $G_2$, then $\alpha_1$ is the short simple root by Convention \ref{conv:G2}.

An element $\delta \in \pi_1( B) = \Hom( \Gm, T)$ is called \emph{minuscule} if $\delta \neq 0$,
and $\langle \alpha, \delta \rangle \in \{0, 1\}$ holds for all positive roots $\alpha \in \Phi^+$ of $G$ with respect to $B$.
Every nonzero element $d \in \pi_1( G)$ has a unique minuscule lift $\delta \in \pi_1( B)$ according to \cite[Chapitre VIII, \S 7, Proposition 8]{bourbaki_Lie_78}.

\begin{lemma} \label{lemma:W_End}
  Suppose that $\delta \in \pi_1( B)$ is zero, or that $-\delta$ is minuscule.
  \begin{itemize}
   \item[i)] If $g \geq 3$, then the complement of $\calW_{\End( \gLie)/k \cdot \id}$ in $\calM_B^{\delta}$ is non-empty.
   \item[ii)] If $g \geq 4$, and $\calW_{\End( \gLie)/k \cdot \id}$ contains a prime divisor $D \subseteq \calM_B^{\delta}$,
    then $\Phi \cong G_2$, $\delta = 0$, and $D = (\alpha_1)_*^{-1}( D')$ for a prime divisor $D' \subseteq \calM_{\Gm}^0$.
  \end{itemize}
\end{lemma}
\begin{proof}
  Note that $\calW_V \subseteq \calW_{V'} \cup \calW_{V''}$ for every short exact sequence
  \begin{equation*}
    0 \longto V' \longto V \longto V'' \longto 0
  \end{equation*}
  of $B$-modules. Using Corollary \ref{cor:End=ext}, we may thus replace the $B$-module $\End(\gLie)/k \cdot \id$
  by the $1$- and $2$-dimensional $B$-modules $V$ listed there.

  The cases $V = \ell_{\alpha}$ and $V = \ell_{\alpha - \beta}$ follow from Lemma \ref{lemma:chi}.
  It remains to treat the cases $V = \wp_i$, $V = \wp^2_i$, and $V = \wp^3_1$.

  The exact sequence \eqref{eq:U'} shows that the kernel of $\pi_i: B \twoheadrightarrow \Ga \rtimes \Gm$ is smooth.
  So the induced map $(\pi_i)_*: \calM_B \to \calM_{\Ga \rtimes \Gm}$ is also smooth, by deformation theory.
  Now use Corollary \ref{cor:b} and Lemma \ref{lemma:Wwp}.
\end{proof}
\begin{proof}[Proof of Proposition \ref{prop:W_End}]
  The claim \ref{prop:W_End}.i follows from Lemma \ref{lemma:W_End}.i.

  Now suppose $g \geq 4$, and that $\calW_{\End( \gLie)/k \cdot \id}$ contains a prime divisor $D \subseteq \calM_G^d$.
  Then its ideal sheaf $\calO( -D)$ is a line bundle, since $\calM_G^d$ is smooth.
  This line bundle is nontrivial because $\rmH^0( \calM_G^d, \calO) = k$, cf. for example \cite[Theorem 5.3.1.i]{pic}.

  For any lift $\delta \in \Hom( \Gm, T)$ of $d$, the pullback of $\calO( -D)$ to $\calM_T^{\delta}$ is still nontrivial,
  due to \cite[Lemma 5.2.6 and Theorem 5.3.1.iv]{pic}. In particular, the pullback of $\calO( -D)$ to $\calM_B^{\delta}$ is also nontrivial.
  Thus the inverse image $\calW_{\End( \gLie)/k \cdot \id} \subseteq \calM_B$ contains a prime divisor in $\calM_B^{\delta}$.

  Assuming moreover that $G$ is of type $G_2$, we have $\pi_1( G) = 0$, so $d = 0$, and $T = \ker( \alpha_1) \times \ker( \alpha_2)$. 
  The pullback of $\calO( -D)$ to $\calM^0_{\ker( \alpha_1)}$ is nontrivial by \cite[Proposition 4.4.7.iii]{pic}.
  Thus $\calW_{\End( \gLie)/k \cdot \id} \subseteq \calM_B$ contains a prime divisor in $\calM_B^0$ that is no pullback along $\alpha_1$.

  These conclusions contradict Lemma \ref{lemma:W_End}. The claim \ref{prop:W_End}.ii follows.
\end{proof}

\section{The obstruction against Poincar\'{e} families} \label{sec:obstructions}
Let the smooth connected group $G$ be reductive, and consider $d \in \pi_1( G)$.
\begin{definition} \label{def:poinc}
  A \emph{Poincar\'{e} family} for an open subscheme $U \subseteq \frakM_G^{d, \regstab}$ is a principal $G$-bundle $\calE$ on $C \times U$ such that
  for every point $x \in U$, the corresponding isomorphism class of stable $G$-bundles contains $\calE|_{C \times \{x\}}$.
\end{definition}
Let $Z \subseteq G$ denote the scheme-theoretic center. The coarse moduli map
\begin{equation*}
  \pi: \calM_G^{d, \regstab} \longto \frakM_G^{d, \regstab}
\end{equation*}
is a gerbe with band $Z$ for the \'{e}tale topology; local sections of $\pi$ are given by the \'{e}tale slices in the proof of Proposition \ref{prop:regstab}. Let
\begin{equation*}
  \psi_G^d \in \Het^2( \frakM_G^{d, \regstab}, Z)
\end{equation*}
be the cohomology class corresponding to the gerbe $\pi$, according to \cite{giraud}.
Note that a Poincar\'{e} family for $U$ is a section of the gerbe $\pi$ over $U$.
Thus $\psi_G^d$ vanishes if and only if there is a Poincar\'{e} family for $\frakM_G^{d, \regstab}$.

Now suppose $\frakM_G^{d, \regstab} \neq \emptyset$.
Slightly abusing notation, we also denote by $\psi_G^d$ its restriction to the generic point $\eta_{\frakM_G^d} \in \frakM_G^{d, \regstab}$. This restriction
\begin{equation*}
  \psi_G^d \in \Het^2( \eta_{\frakM_G^d}, Z)
\end{equation*}
vanishes if and only if there is a Poincar\'{e} family for some open subscheme $\emptyset \neq U \subseteq \frakM_G^{d, \regstab}$.
In this section, we determine the order of these obstruction classes. The result is given in terms of the root system of $G$.

Let $G$ for the moment be of adjoint type. Choose a maximal torus $T \subseteq G$. Let $\Lambda_{\coroots} \subseteq \Lambda_T := \Hom( \Gm, T)$ denote
the subgroup generated by the coroots of $G$. The Weyl group $W$ of $(G, T)$ acts on $\Lambda_T$.
This action preserves the subgroup $\Lambda_{\coroots}$, and the induced action on the quotient group $\Lambda_T/\Lambda_{\coroots} = \pi_1( G)$ is trivial.
\begin{lemma}
  Given an even $W$-invariant symmetric bilinear form
  \begin{equation*}
    \Lambda_{\coroots} \times \Lambda_{\coroots} \longto \Z,
  \end{equation*}
  its bilinear extension $\Lambda_T \times \Lambda_T \to \Q$ descends to a symmetric bilinear map
  \begin{equation*}
    \pi_1( G) \times \pi_1( G) \longto \Q/\Z.
  \end{equation*}
\end{lemma}
\begin{proof}
  Let $b: \Lambda_{\coroots} \times \Lambda_{\coroots} \to \Z$ be $W$-invariant, symmetric, bilinear,
  and even in the sense that $b( \lambda, \lambda)$ is even for all $\lambda \in \Lambda_{\coroots}$.
  Given a root $\alpha: \Lambda_T \to \Z$ of $G$, with corresponding coroot $\alpha^{\vee} \in \Lambda_T$, the formula
  \begin{equation*}
    b( \lambda, \alpha^{\vee}) = \alpha( \lambda) \cdot b( \alpha^{\vee}, \alpha^{\vee})/2
  \end{equation*}
  holds for all $\lambda \in \Lambda_{\coroots}$, according to \cite[Chapitre VI, \S 1, Lemme 2]{bourbaki_lie_456}.
  In particular, the map $b(\_,\alpha^{\vee}): \Lambda_{\coroots} \to \Z$ is an integral multiple of $\alpha$, and hence extends to a linear map $\Lambda_T \to \Z$.
  Thus the bilinear extension $\Lambda_T \times \Lambda_T \to \Q$ of $b$ is integral on $\Lambda_T \times \Lambda_{\coroots}$,
  and also on $\Lambda_{\coroots} \times \Lambda_T$ by symmetry. Hence the composition $\Lambda_T \times \Lambda_T \to \Q \twoheadrightarrow \Q/\Z$
  descends to a bilinear map $\pi_1( G) \times \pi_1( G) \to \Q/\Z$.
\end{proof}
\begin{definition}
  Suppose that $G$ is of adjoint type. We denote by
  \begin{equation*}
    \Psi( G) \subseteq \Hom( \pi_1( G) \otimes \pi_1( G), \Q/Z)
  \end{equation*}
  the abelian group of all bilinear maps $b: \pi_1( G) \times \pi_1( G) \to \Q/\Z$
  that come from even $W$-invariant symmetric bilinear forms $\Lambda_{\coroots} \times \Lambda_{\coroots} \to \Z$.
\end{definition}
Note that $\Psi( G)$ is determined by the root system of $G$. If $G = G_1 \times G_2$, then $\Psi( G) = \Psi( G_1) \oplus \Psi( G_2)$.
If $G$ is a simple group, then the abelian group $\Psi( G)$ is cyclic, and a generator is given by Table \ref{tab:Psi}.
\renewcommand{\arraystretch}{1.5}
\begin{table} \begin{tabular}{c|c|c|cc} 
  type of $G$ & $\pi_1( G)$ & $\Psi( G)$ & generator of $\Psi( G)$ & \\\hline \hline
  $A_l$, $l \geq 1$ & $\mods{l+1}$ & $\mods{l+1}$ & mult: $\mods{l+1} \otimes \mods{l+1} \longto \mods{l+1}$ &\\\hline
  $B_l$, $l \geq 2$ & $\mods{2}$ & $0$ & --- &\\\hline
  & & $\mods{2}$ & mult: $\mods{2} \otimes \mods{2} \longto \mods{2}$ & $l$ odd\\
  \raisebox{2ex}[-2ex]{$C_l$, $l \geq 3$} & \raisebox{2ex}[-2ex]{$\mods{2}$} & $0$ & --- & $l$ even\\\hline
  & $\mods{4}$ & $\mods{4}$ & mult: $\mods{4} \otimes \mods{4} \longto \mods{4}$ & $l$ odd\\
  $D_l$, $l \geq 4$ && & $(\begin{smallmatrix} 0&1\\1&0 \end{smallmatrix}): (\mods{2})^2 \otimes (\mods{2})^2 \longto \mods{2}$ & $l \in 4 \Z$\\
  & \raisebox{2ex}[-2ex]{$(\mods{2})^2$} & \raisebox{2ex}[-2ex]{$\mods{2}$}
    & $(\begin{smallmatrix} 1&0\\0&1 \end{smallmatrix}): (\mods{2})^2 \otimes (\mods{2})^2 \longto \mods{2}$ & $l \in 4\Z + 2$\\ \hline
  $E_6$ & $\mods{3}$ & $\mods{3}$ & mult: $\mods{3} \otimes \mods{3} \longto \mods{3}$ &\\\hline
  $E_7$ & $\mods{2}$ & $\mods{2}$ & mult: $\mods{2} \otimes \mods{2} \longto \mods{2}$ &\\\hline
  $E_8$ & $0$ & $0$ & --- &\\\hline
  $F_4$ & $0$ & $0$ & --- &\\\hline
  $G_2$ & $0$ & $0$ & --- &
\end{tabular} \caption[The abelian group $\Psi( G)$ for simple groups $G$]{\label{tab:Psi} The abelian group $\Psi( G)$ for simple $G$ of adjoint type.\newline
(To obtain the required maps to $\Q/\Z$, embed $\Z/n$ into $\Q/\Z$).}
\end{table}

Now return to the general case where $G$ is reductive. Let $Z^0$ be the (reduced) identity component in the center $Z$ of $G$. The central isogenies
\begin{equation*}
  G' := [G, G] \twoheadrightarrow \Gbar := G/Z^0 \twoheadrightarrow \Gad := G/Z
\end{equation*}
correspond to subgroups $\pi_1( G') \subseteq \pi_1( \Gbar) \subseteq \pi_1( \Gad)$.
If $G$ is semisimple, then $G' = \Gbar = G$, so we just have one subgroup $\pi_1( G) \subseteq \pi_1( \Gad)$.
\begin{definition} \label{def:Psi'}
  For reductive $G$, we denote by $\Psi'( G) \subseteq \Psi( \Gad)$ the subgroup
  of all elements $b: \pi_1( \Gad) \times \pi_1( \Gad) \to \Q/\Z$ in $\Psi( \Gad)$ with
  \begin{equation*}
    b( \pi_1( \Gbar) \times \pi_1( G')) = 0.
  \end{equation*}
  Given an element $d \in \pi_1( G)$ with image $\dbar \in \pi_1( \Gbar)$, we denote by
  \begin{equation*}
    \ev_G^d: \Psi'( G) \longto \Hom( \frac{\pi_1( \Gad)}{\pi_1( G')}, \Q/\Z)
  \end{equation*}
  the evaluation map that sends $b$ to $b( \dbar, \_): \pi_1( \Gad) \to \Q/Z$.
\end{definition}
\begin{remark}
  This finite abelian group $\Hom( \frac{\pi_1( \Gad)}{\pi_1( G')}, \Q/\Z)$ can be identified with the character group $\Hom(Z', \Gm)$ of the center
  $Z' \subseteq G'$; cf. the equations \eqref{eq:Z'} in the proof of Proposition \ref{prop:weights} below.
\end{remark}
\begin{theorem} \label{thm:ob_generic}
  Assume $g \geq 3$. Let $G$ be reductive, with scheme-theoretic center $Z$.
  The order of $\psi_G^d$ in $\Het^2( \eta_{\frakM_G^d}, Z)$ is the exponent of the finite abelian group $\coker( \ev_G^d)$ in Definition \ref{def:Psi'}.
\end{theorem}
\begin{corollary} \label{cor:ob_generic}
  Assume $g \geq 3$. There is a Poincar\'{e} family for some non-empty open subscheme $U \subseteq \frakM_G^d$
  if and only if the evaluation map $\ev_G^d$ in Definition \ref{def:Psi'} is surjective.
\end{corollary}
\begin{theorem} \label{thm:ob_global}
  Assume $g \geq 4$. The order of $\psi_G^d$ in $\Het^2( \frakM_G^{d, \regstab}, Z)$ is the least common multiple of its order in $\Het^2( \eta_{\frakM_G^d}, Z)$
  and the exponent of the finite abelian group $\Hom( Z/Z^0, \Gm)$.
\end{theorem}
\begin{corollary} \label{cor:ob_global}
  Assume $g \geq 4$, and that there is a Poincar\'{e} family for some non-empty open subscheme $U \subseteq \frakM_G^d$.
  Then there is a Poincar\'{e} family for $\frakM_G^{d, \regstab}$ if and only if $Z$ is a torus.
\end{corollary}
Theorem \ref{thm:ob_generic} and Theorem \ref{thm:ob_global} are proved in the next section.
\begin{remark}
  The assumption $g \geq 3$ (respectively, $g \geq 4$) is needed only to ensure that $\calM_G^{d, \regstab}$ is non-empty
  (respectively, that its complement has codimension $\geq 2$ in $\calM_G^d$).
  In the case $\characteristic( k) = 0$, these statements about $\calM_G^{d, \regstab}$ are proved also for smaller $g$ in \cite[Theorem II.6]{faltings}.
  In this case, Theorem \ref{thm:ob_generic} and its Corollary \ref{cor:ob_generic} follow also for $g = 2$, while
  Theorem \ref{thm:ob_global} and its Corollary \ref{cor:ob_global} follow also for $g = 3$,
  and even for $g = 2$ unless there is a nontrivial homomorphism $G \to \PGL_2$.
\end{remark}

\section{Gerbes with band $Z$} \label{sec:gerbes}

We say that an algebraic group $Z$ over $k$ is of multiplicative type if it is of the form $Z \cong \Gm^r \times \mu_{n_1} \times \cdots \times \mu_{n_s}$
with $r, s \geq 0$ and $n_1, \ldots, n_s \geq 1$.

Let $U$ be an integral scheme of finite type over $k$. Let $\pi: \calU \to U$ be a gerbe with band $Z$ for an algebraic group $Z$ of multiplicative type over $k$.
We denote the cohomology class of this gerbe by $\psi_{\calU} \in \Het^2( U, Z)$.

The stack $\calU$ is given by a groupoid $\calU( S)$ for each $k$-scheme $S$.
Since $\pi$ is a gerbe with band $Z$, we have an isomorphism $\iota_{\calE}: Z( S) \to \Aut_{\calU( S)}( \calE)$ for every object $\calE$ in $\calU( S)$.
These data are subject to appropriate compatibility conditions.

Now let $\calL$ be a line bundle on $\calU$. Recall that $\calL$ is given by a functor $\calL_S$ from $\calU( S)$ to the groupoid
of line bundles on $S$ for every $k$-scheme $S$. In particular, $\calL_S$ defines for every object $\calE$ in $\calU( S)$ a group homomorphism
\begin{equation*}
  \calL_{S, \calE}: \Aut_{\calU( S)}( \calE) \longto \Aut_{\calO_S}( \calL_S( \calE)) = \Gamma( S, \calO_S^*).
\end{equation*}
The compatibility conditions ensure that the compositions
\begin{equation*}
  Z( S) \longto[ \iota_{\calE}] \Aut_{\calU( S)}( \calE) \longto[ \calL_{S, \calE}] \Gamma( S, \calO_S^*)
\end{equation*}
define a $1$-morphism $Z \times \calU \to \Gm \times \calU$ over $\calU$.
Because $\calU$ is connected and $\Hom( Z, \Gm)$ is discrete, this $1$-morphism is the pullback of some character $\chi: Z \to \Gm$.
We call $\chi$ the weight of $\calL$. Sending each line bundle $\calL$ on $\calU$ to its weight $\chi$ defines a group homomorphism
\begin{equation*}
  \weight: \Pic( \calU) \longto \Hom( Z, \Gm).
\end{equation*}
If $\calL$ has trivial weight, then it descend to $U$, so we have an exact sequence
\begin{equation*}
  0 \longto \Pic( U) \longto \Pic( \calU) \longto[\weight] \Hom( Z, \Gm).
\end{equation*}
\begin{example}
  Let $\emptyset \neq \calU \subseteq \calM_G^d$ be an open substack, with $G$ reductive.
  Given a $k$-scheme $S$, the objects in $\calU( S)$ are principal $G$-bundles $\calE$ on $C \times S$.
  The center $Z \subseteq G$ acts by automorphisms on every principal $G$-bundle.
  Thus we obtain a homomorphism $\iota_{\calE}: Z( S) \to \Aut_{\calU( S)}( \calE)$ for every object $\calE$ in $\calU( S)$.

  Suppose $\calU \subseteq \calM_G^{d, \regstab}$, and let $U \subseteq \frakM_G^{d, \regstab}$ be the corresponding open subscheme.
  Then all $\iota_{\calE}$ are isomorphisms, and they turn the coarse moduli map $\pi: \calU \to U$ into a gerbe with band $Z$.
  In particular, the action of $Z$ on fibers of line bundles $\calL$ on $\calU$ defines a group homomorphism
  \begin{equation*}
    \weight = \weight_G^d: \Pic( \calU) \longto \Hom( Z, \Gm).
  \end{equation*}
  However, the definition of $\weight( \calL)$ used only the homomorphisms $\iota_{\calE}$, not the fact that they are isomorphisms.
  Thus we can drop the assumption $\calU \subseteq \calM_G^{d, \regstab}$, and the same construction yields a group homomorphism
  \begin{equation*}
    \weight = \weight_G^d: \Pic( \calU) \longto \Hom( Z, \Gm)
  \end{equation*}
  for every non-empty open substack $\calU \subseteq \calM_G^d$.
\end{example}
\begin{proposition} \label{prop:weights}
  Given $G$ reductive and $d \in \pi_1( G)$, consider the map
  \begin{equation*}
    \weight = \weight_G^d: \Pic( \calM_G^d) \longto \Hom( Z, \Gm).
  \end{equation*}
  \begin{itemize}
   \item[i)] The kernel of $\weight_G^d$ contains the torsion in $\Pic( \calM_G^d)$.
   \item[ii)] The cokernel of $\weight_G^d$ is isomorphic to the cokernel of the map $\ev_G^d$ in Definition \ref{def:Psi'}.
  \end{itemize}
\end{proposition}
\begin{proof}
  In \cite[Definition 5.2.1]{pic}, we have defined a finitely generated free abelian group $\NS( \calM_G^d)$ in terms of the root system of $G$.
  Instead of repeating the definition, we just recall the properties that we need here:
  \begin{itemize}
   \item[(a)] There is a canonical epimorphism $c_G: \Pic( \calM_G^d) \twoheadrightarrow \NS( \calM_G^d)$.
   \item[(b)] Every homomorphism of reductive groups $\varphi: G \to H$ induces a homomorphism of abelian groups $\varphi^{\NS, d}$
    such that the diagram
    \begin{equation*} \xymatrix{
      \Pic( \calM_H^e) \ar@{->>}[r]^{c_H} \ar[d]_{\varphi^*} & \NS( \calM_H^e) \ar[d]^{\varphi^{\NS, d}}\\
      \Pic( \calM_G^d) \ar@{->>}[r]^{c_G} & \NS( \calM_G^d)
    } \end{equation*}
    commutes, with $e := \varphi_*( d) \in \pi_1( H)$.
   \item[(c)] If $G = T$ is a torus, then $\NS( \calM_T^d) = \Hom( T, \Gm) \oplus \NS( \frakM_T^0)$, and the weight map
    $\weight_T^d: \Pic( \calM_T^d) \to \Hom( T, \Gm)$ is the composition
    \begin{equation*}
      \Pic( \calM_T^d) \xrightarrow{c_T} \NS( \calM_T^d) \xrightarrow{\pr_1} \Hom( T, \Gm).
    \end{equation*}
   \item[(d)] Choose a maximal torus $T \subseteq G$. In the notation of Definition \ref{def:Psi'}, let $T' \twoheadrightarrow \bar{T} \twoheadrightarrow T^{\ad}$
    be the induced tori in $G' \twoheadrightarrow \bar{G} \twoheadrightarrow \Gad$. Their lattices
    \begin{equation*}
      \Lambda_{T'} := \Hom( \Gm, T') \subseteq \Lambda_{\bar{T}} := \Hom( \Gm, \bar{T}) \subseteq \Lambda_{T^{\ad}} := \Hom( \Gm, T^{\ad})
    \end{equation*}
    all contain $\Lambda_{\coroots}$. Let $\widetilde{\Psi}( G)$ denote the group of all bilinear maps
    \begin{equation*}
      b: \Lambda_{\bar{T}} \times \Lambda_{T'} \longto \Z
    \end{equation*}
    whose restriction $\Lambda_{\coroots} \times \Lambda_{\coroots} \to \Z$ is $W$-invariant, symmetric, and even.
    Then one has a canonical exact sequence
    \begin{equation*}
      0 \longto \NS( \calM_{G/G'}^e) \xrightarrow{\pi^{\NS, d}} \NS( \calM_G^d) \xrightarrow{\pr_2} \widetilde{\Psi}( G) \longto 0
    \end{equation*}
    where $\pi: G \twoheadrightarrow G/G'$ is the projection, and $e := \pi_*( d) \in \pi_1( G/G')$.
   \item[(e)] Let $\iota: T \hookrightarrow G$ be the inclusion. Choose a lift $\delta \in \Lambda_T$ of $d \in \pi_1( G)$, and let
    ${\bar{\delta}} \in \Lambda_{\bar{T}}$ be the image of $\delta$. Then the following diagram commutes:
    \begin{equation*} \xymatrix{
      \NS( \calM_G^d) \ar[rr]^{\pr_2} \ar[d]_{\pr_1 \circ \iota^{\NS, \delta}} && \widetilde{\Psi}( G) \ar[d]^{b \mapsto b( -\bar{\delta}, \_)}\\
      \Hom( T, \Gm) \ar@{->>}[rr]^{\chi \mapsto \chi|_{T'}} && \Hom( T', \Gm)
    } \end{equation*}
  \end{itemize}
  These properties are contained in Proposition 5.2.11 and Theorem 5.3.1 of \cite{pic}; see also its Subsection 3.2 for (c),
  and its Definition 5.2.5 for (e).

  The inclusion $\iota: T \hookrightarrow G$ of a maximal torus and the maximal commutative quotient $\pi: G \twoheadrightarrow G/G'$ induce a commutative diagram
  \begin{equation} \label{diag:weights} \xymatrix{
    \Pic(\calM_{G/G'}^e) \ar[r]^{\pi^*} \ar[d]^{\weight_{G/G'}^e} & \Pic(\calM_G^d) \ar[r]^{\iota^*} \ar[d]^{\weight_G^d} & \Pic(\calM_T^{\delta}) \ar[d]^{\weight_T^{\delta}}\\
    \Hom( G/G', \Gm) \ar[r]^-{\pi^*} & \Hom( Z, \Gm) & \Hom( T, \Gm). \ar@{->>}[l]
  } \end{equation}
  In particular, the map $\weight_G^d$ factors through the torsionfree abelian group $\Hom( T, \Gm)$; this proves part (i).

  Let $Z'$ denote the scheme-theoretic center of $G'$; then $Z/Z' \cong G/G'$. In this diagram \eqref{diag:weights},
  the map $\weight_{G/G'}^e$ is surjective by property (c), so the image of $\weight_G^d$ contains the image $\Hom( Z/Z', \Gm)$ of $\pi^*$;
  consequently, the cokernel of $\weight_G^d$ is isomorphic to the cokernel of the composition
  \begin{equation} \label{eq:composition}
    \Pic( \calM_G^d) \xrightarrow{\iota^*} \Pic( \calM_T^{\delta}) \xrightarrow{\weight_T^{\delta}} \Hom(T, \Gm) \twoheadrightarrow \Hom(Z', \Gm).
  \end{equation}
  Using the canonical identification
  \begin{equation} \label{eq:Z'}
    \Hom( Z', \Gm) = \frac{\Hom( \Lambda_{T'}, \Z)}{\Hom( \Lambda_{T^{\ad}}, \Z)}
      = \Hom( \frac{\Lambda_{T^{\ad}}}{\Lambda_{T'}}, \frac{\Q}{\Z}) = \Hom( \frac{\pi_1( \Gad)}{\pi_1( G')}, \frac{\Q}{\Z}),
  \end{equation}
  the above properties (a) -- (e) yield a commutative diagram
  \begin{equation*} \xymatrix{
    \Pic( \calM_G^d) \ar@{->>}[r]^{c_G} \ar[d]^{\iota^*} & \NS( \calM_G^d) \ar@{->>}[r]^{\pr_2} \ar[d]^{\pr_1 \circ \iota^{\NS, \delta}}
      & \widetilde{\Psi}( G) \ar@{->>}[r] \ar[d]^{b \mapsto b( -\bar{\delta}, \_)} & \Psi'( G) \ar[d]^{-\ev_G^d}\\
    \Pic( \calM_T^{\delta}) \ar[r]^-{\weight_T^{\delta}} & \Hom( T, \Gm) \ar@{->>}[r] & \Hom( T', \Gm) \ar@{->>}[r] & \Hom( Z', \Gm)
  } \end{equation*}
  where the unlabeled arrows are the obvious canonical surjections. Because the three maps in the top row of this diagram are surjective,
  the cokernel of the composition \eqref{eq:composition} is isomorphic to the cokernel of $\ev_G^d$.
\end{proof}
\begin{samepage} \begin{lemma} \label{lemma:extend_L}
  Let $\calM$ be an Artin stack which is locally of finite type and smooth over $k$. Let $\calU \subseteq \calM$ be an open substack.
  \begin{itemize}
   \item[i)] The restriction map $\Pic( \calM) \to \Pic( \calU)$ is surjective.
   \item[ii)] If $\calM \setminus \calU$ has codimension $\geq 2$ in $\calM$, then $\Pic( \calM) = \Pic( \calU)$.
  \end{itemize}
\end{lemma} \end{samepage}
\begin{proof}
  i) Suppose $\calU \subsetneq \calM$, and let $\calL$ be a line bundle on $\calU$.
  Using Zorn's lemma, it suffices to extend $\calL$ to a line bundle on $\calM'$ for some open substack $\calU \subsetneq \calM' \subseteq \calM$.
  We can choose a quasi-compact open substack $\calU' \subseteq \calM$ not contained in $\calU$, and take $\calM' := \calU \cup \calU'$.
  To extend $\calL$ to $\calM'$, it suffices to extend $\calL|_{\calU \cap \calU'}$ to $\calU'$.
  Thus we may assume without loss of generality that $\calM$ is quasi-compact.

  Now $\calM$ and $\calU$ are of finite type over $k$. Then \cite[Corollaire 15.5]{L-MB} allows us to extend $\calL$ to a coherent sheaf on $\calM$,
  namely to a coherent subsheaf of $j_*( \calL)$, where $j: \calU \hookrightarrow \calM$ is the inclusion.
  Using smoothness, the bidual of this coherent sheaf of rank one on $\calM$ is the required line bundle extending $\calL$.

  ii) Suppose that $\calM \setminus \calU$ has codimension $\geq 2$ in $\calM$, and let $\calL$ be a line bundle on $\calM$.
  Every section of $\calL$ or $\calL^{-1}$ over $\calU$ extends uniquely to $\calM$ by Hartog's theorem.
  Hence the restriction map $\Pic( \calM) \to \Pic( \calU)$ is also injective.
\end{proof}
\begin{lemma} \label{lemma:annihilated}
  Let $\varphi: Z \to Z'$ be a homomorphism of algebraic groups over $k$, with $Z$ and $Z'$ of multiplicative type.
  Given an integral scheme $U$ of finite type over $k$, and a gerbe $\pi: \calU \to U$ with band $Z$, the following conditions are equivalent:
  \begin{itemize}
   \item[i)] The class $\psi_{\calU} \in \Het^2( U, Z)$ is in the kernel of $\Het^2( U, Z) \longto[ \varphi_*] \Het^2( U, Z')$.
   \item[ii)] There is a group homomorphism $\sigma: \Hom( Z', \Gm) \to \Pic( \calU)$ with
    \begin{equation*}
      \weight \circ \sigma = \varphi^*: \Hom( Z', \Gm) \longto \Hom( Z, \Gm).
    \end{equation*}
  \end{itemize}
\end{lemma}
\begin{proof}
  Let $\pi': \calU' \to U$ be a gerbe with band $Z'$, and let $\psi_{\calU'} \in \Het^2( U, Z')$ be its class.
  According to Proposition IV.3.1.5 and Th\'{e}or\`{e}me IV.3.4.2 in \cite{giraud}, the relation $\varphi_*( \psi_{\calU}) = \psi_{\calU'}$ holds
  if and only if there is a $1$-morphism $\calU \to \calU'$ over $U$ which induces $\varphi$ on automorphism groups.

  In particular, $\psi_{\calU}$ is in the kernel of $\varphi_*$ if and only if
  there is a $1$-morphism $\Sigma: \calU \to BZ'$ which induces $\varphi$ on automorphism groups.

  Given such a $1$-morphism $\Sigma$, the required group homomorphism $\sigma$ is
  \begin{equation*}
    \sigma := \Sigma^*: \Hom( Z', \Gm) = \Pic( BZ') \longto \Pic( \calU).
  \end{equation*}

  Conversely, let $\sigma: \Hom( Z', \Gm) \to \Pic( \calU)$ be given. Choosing an isomorphism
  \begin{equation*}
    Z' \cong \Gm^r \times \mu_{n_1} \times \cdots \times \mu_{n_s},
  \end{equation*}
  let $\chi_1, \ldots, \chi_r: Z' \to \Gm$ and $\chi_{r+i}: Z' \to \mu_{n_i}$ for $i = 1, \ldots, s$ be the projections.
  Their images under $\sigma$ are line bundles $\calL_1, \ldots, \calL_{r+s}$ on $\calU$ such that $\calL_{r+i}^{\otimes n_i}$ is trivial for all $i$.

  The tuple $(\calL_1, \ldots, \calL_r)$ defines a $1$-morphism $\calU \to B\Gm^r$, and $\calL_{r+i}$ together with a trivialisation of $\calL_{r+i}^{\otimes n_i}$
  defines a $1$-morphism $\calU \to B\mu_{n_i}$. Let $\Sigma: \calU \to BZ'$ be the product of these $1$-morphisms.
  If $\weight \circ \sigma = \varphi^*$, then $\Sigma$ induces $\varphi$ on automorphism groups, and $\varphi_*( \psi_{\calU}) = 0$ follows.
\end{proof}
\begin{proof}[Proof of Theorem \ref{thm:ob_generic}]
  Consider a positive integer $n \geq 1$.

  Suppose $n \cdot \psi_G^d = 0$ generically. Then $n \cdot \psi_{\calU} = 0 \in \Het^2( U, Z)$ for some open subscheme
  $\emptyset \neq U \subseteq \frakM_G^{d, \regstab}$, with $\calU := \pi^{-1}( U) \subseteq \calM_G^{d, \regstab}$.
  Applying Lemma \ref{lemma:annihilated} to $n: Z \to Z$, we get a map
  \begin{equation*}
    \sigma: \Hom(Z, \Gm) \longto \Pic( \calU)
  \end{equation*}
  with $\weight_G^d \circ \sigma = n$. In particular, the cokernel of
  \begin{equation*}
    \weight_G^d: \Pic( \calU) \longto \Hom( Z, \Gm)
  \end{equation*}
  is annihilated by $n$. But $\weight_G^d: \Pic( \calM_G^d) \to \Hom( Z, \Gm)$ has the same cokernel, because the restriction map
  $\Pic( \calM_G^d) \to \Pic( \calU)$ is surjective due to Lemma \ref{lemma:extend_L}.i.
  Using Proposition \ref{prop:weights}.ii, it follows that the cokernel of $\ev_G^d$ is also annihilated by $n$.

  Conversely, suppose that $\coker( \ev_G^d)$ is annihilated by $n$. We write
  \begin{equation*}
    \Hom(Z, \Gm) \cong \Z^r \oplus \Z/n_1 \oplus \cdots \oplus \Z/n_s,
  \end{equation*}
  and let $\chi_1, \ldots, \chi_{r+s} \in \Hom(Z, \Gm)$ be the corresponding generators.
  Using the assumption, Proposition \ref{prop:weights}.ii allows us to find line bundles
  \begin{equation*}
    \calL_1, \ldots, \calL_{r+s} \in \Pic( \calM_G^d)
  \end{equation*}
  such that $\calL_i$ has weight $n \chi_i$. Then the line bundles
  \begin{equation*}
    \calL_{r+1}^{\otimes n_1}, \ldots, \calL_{r+s}^{\otimes n_s}
  \end{equation*}
  all have trivial weight. Hence their restrictions descend to $\frakM_G^{d, \regstab}$, and are thus trivial over some open subscheme
  $\emptyset \neq U \subseteq \frakM_G^{d, \regstab}$.
  Consequently, the line bundles $\calL_1, \ldots, \calL_{r+s}$ define a group homomorphism
  \begin{equation*}
    \sigma: \Hom( Z, \Gm) \longto \Pic( \calU), \qquad \calU := \pi^{-1}( U) \subseteq \calM_G^{d, \regstab}.
  \end{equation*}
  We have $\weight_G^d \circ \sigma = n$ by construction. Hence $n \cdot \psi_{\calU} = 0 \in \Het^2( U, Z)$ according to Lemma \ref{lemma:annihilated}.
  This shows $n \cdot \psi_G^d = 0$ generically.
\end{proof}

\begin{proof}[Proof of Theorem \ref{thm:ob_global}]
  Consider a positive integer $n \geq 1$.

  Suppose $n \cdot \psi_G^d = 0$ in $\Het^2( \frakM_G^{d, \regstab}, Z)$. Applying Lemma \ref{lemma:annihilated} to the map $n: Z \to Z$,
  we get a homomorphism $\sigma: \Hom(Z, \Gm) \to \Pic( \calM_G^{d, \regstab})$ with $\weight_G^d \circ \sigma = n$.
  Here the map $\weight_G^d$ vanishes on all torsion elements, according to Proposition \ref{prop:weights}.i and Lemma \ref{lemma:extend_L}.ii.
  It follows that $n$ annihilates the torsion subgroup $\Hom( Z/Z^0, \Gm)$ of $\Hom( Z, \Gm)$.

  Conversely, suppose that $n$ annihilates $\Hom( Z/Z^0, \Gm)$, and $n \cdot \psi_G^d = 0$ generically.
  Then there is an open subscheme $\emptyset \neq U \subseteq \frakM_G^{d, \regstab}$ with
  \begin{equation*}
    n \cdot \psi_{\calU} = 0 \in \Het^2( U, Z), \qquad \calU := \pi^{-1}( U) \subseteq \calM_G^{d, \regstab}.
  \end{equation*}
  We choose characters $\chi_1, \ldots, \chi_r: Z \to \Gm$ whose restrictions to $Z^0$ form an isomorphism $Z^0 \cong \Gm^r$.
  Applying Lemma \ref{lemma:annihilated} to $n: Z \to Z$, and using Lemma \ref{lemma:extend_L}.i, we can find line bundles
  \begin{equation*}
    \calL_1, \ldots, \calL_r \in \Pic( \calM_G^{d, \regstab})
  \end{equation*}
  such that $\calL_i$ has weight $n \chi_i$. We define $\sigma^0: \Hom( Z^0, \Gm) \to \Pic( \calM_G^{d, \regstab})$ by sending $\chi_i$ to $\calL_i$.
  Then the homomorphism
  \begin{equation*}
    \sigma: \Hom( Z, \Gm) \longto \Pic( \calM_G^{d, \regstab}), \qquad \chi \mapsto \sigma^0( \chi|_{Z^0}),
  \end{equation*}
  satisfies the equation $\weight_G^d \circ \sigma = n$, because both sides of this equation vanish on $\Hom( Z/Z^0, \Gm)$ and map $\chi_i$ to $n \chi_i$.

  Due to Lemma \ref{lemma:annihilated}, this implies $n \cdot \psi_G^d = 0$ in $\Het^2( \frakM_G^{d, \regstab}, Z)$.
\end{proof}


\begin{thebibliography}{10}

\bibitem{BBNN}
V.~Balaji, I.~Biswas, D.S. Nagaraj, and P.E. Newstead.
\newblock {Universal families on moduli spaces of principal bundles on curves.}
\newblock {\em Int. Math. Res. Not.}, Article ID 80641, 2006.

\bibitem{bardsley-richardson}
P.~Bardsley and R.W. Richardson.
\newblock {\'Etale slices for algebraic transformation groups in characteristic
  p.}
\newblock {\em Proc. Lond. Math. Soc., III. Ser.}, 51:295--317, 1985.

\bibitem{pic}
I.~Biswas and N.~Hoffmann.
\newblock {The line bundles on moduli stacks of principal bundles on a curve.}
\newblock {\em Documenta Math.}, 15:35--72, 2010.

\bibitem{borel}
A.~Borel.
\newblock {Properties and linear representations of Chevalley groups.}
\newblock In {\em {Seminar Alg. Groups and Related Finite Groups, Princeton
  1968/69}}, {Springer Lecture Notes in Mathematics 131}, 1970.

\bibitem{bourbaki_lie_456}
N.~Bourbaki.
\newblock {\em {\'El\'ements de math\'ematique. Groupes et alg\`ebres de Lie.
  Chapitres IV, V et VI.}}
\newblock {Paris: Hermann}, 1968.

\bibitem{bourbaki_Lie_78}
N.~Bourbaki.
\newblock {\em {\'El\'ements de math\'ematique. Groupes et alg\`ebres de Lie.
  Chapitres VII et VIII.}}
\newblock {Paris: Hermann}, 1975.

\bibitem{carter}
R.W. Carter.
\newblock {\em {Simple groups of Lie type.}}
\newblock {London etc.: John Wiley \& Sons}, 1972.

\bibitem{demazure-gabriel}
M.~Demazure and P.~Gabriel.
\newblock {\em {Groupes alg\'ebriques.}}
\newblock {Paris: Masson et Cie, \'Editeur; Amsterdam: North-Holland Publishing
  Company}, 1970.

\bibitem{sga3.3}
M.~Demazure and A.~Grothendieck.
\newblock {\em {SGA 3: Sch\'emas en groupes. Expos\'es XIX \`a XXVI.}}
\newblock {Springer Lecture Notes in Mathematics 153}, 1970.

\bibitem{faltings}
G.~Faltings.
\newblock {Stable $G$-bundles and projective connections.}
\newblock {\em J. Algebr. Geom.}, 2(3):507--568, 1993.

\bibitem{giraud}
J.~Giraud.
\newblock {\em {Cohomologie non ab\'{e}lienne.}}
\newblock {Grundlehren der mathematischen Wissenschaften, Band 179.
  Berlin-Heidelberg-New York: Springer-Verlag}, 1971.

\bibitem{GLSS}
T.L. G\'{o}mez, A.~Langer, A.H.W. Schmitt, and I.~Sols.
\newblock Moduli spaces for principal bundles in arbitrary characteristic.
\newblock {\em Adv. Math.}, 219(4):1177--1245, 2008.

\bibitem{GLSS2}
T.L. G\'{o}mez, A.~Langer, A.H.W. Schmitt, and I.~Sols.
\newblock Moduli spaces for principal bundles in large characteristic.
\newblock In I.~Biswas, R.S. Kulkarni, and S.~Mitra, editors, {\em Proceedings
  of `Teichm{\"u}ller Theory and Moduli Problems' (Allahabad, 2006)}, pages
  281--371, 2010.
\newblock Ramanujan Mathematical Society Lecture Notes Series, Vol. 10.

\bibitem{jochen}
J.~Heinloth.
\newblock {Semistable reduction for $G$-bundles on curves.}
\newblock {\em J. Algebr. Geom.}, 17:167--183, 2008.

\bibitem{jochen2}
J.~Heinloth.
\newblock {Addendum to "Semistable reduction for $G$-bundles on curves."}.
\newblock {\em J. Algebr. Geom.}, 19:193--197, 2010.

\bibitem{par}
N.~Hoffmann.
\newblock Rationality and poincar\'{e} families for vector bundles with extra
  structure on a curve.
\newblock {\em Int. Math. Res. Not.}, 2007.
\newblock article ID rnm010, 29 pages.

\bibitem{pi0}
N.~Hoffmann.
\newblock {On Moduli Stacks of $G$-Bundles over a Curve}.
\newblock to appear in: A.H.W. Schmitt (Ed.), Proceedings of `Affine Flag
  Manifolds and Principal Bundles' (Berlin 2008). {Birkh\"auser} Trends in
  Mathematics, 2010.
\newblock preprint available at \verb|http://userpage.fu-berlin.de/~nhoffman|.

\bibitem{yogish}
Y.~Holla.
\newblock Parabolic reductions of principal bundles.
\newblock preprint arXiv:math/0204219.
\newblock available at \verb|http://www.arXiv.org|.

\bibitem{L-MB}
G.~Laumon and L.~Moret-Bailly.
\newblock {\em {Champs alg\'ebriques.}}
\newblock {Ergebnisse der Mathematik und ihrer Grenzgebiete, Band 39. Berlin:
  Springer}, 2000.

\bibitem{luna}
D.~Luna.
\newblock {Slices \'etal\'es.}
\newblock {\em Bull. Soc. Math. Fr., Suppl., Mem.}, 33:81--105, 1973.

\bibitem{Ne}
P.E. Newstead.
\newblock {A nonexistence theorem for families of stable bundles.}
\newblock {\em J. Lond. Math. Soc., II. Ser.}, 6:259--266, 1973.

\bibitem{Ra}
S.~Ramanan.
\newblock {The moduli space of vector bundles over an algebraic curve.}
\newblock {\em Math. Ann.}, 200:69--84, 1973.

\bibitem{raynaud}
M.~Raynaud.
\newblock {Sections des fibr\'{e}s vectoriels sur une courbe.}
\newblock {\em Bull. Soc. Math. Fr.}, 110:103--125, 1982.

\end{thebibliography}
\end{document}